\documentclass{ifacconf}

\usepackage{graphicx}      
\usepackage{natbib}        

\usepackage{amsmath}
\usepackage{amssymb}
\usepackage{mathtools}

\newtheorem{definition}{Definition}
\newtheorem{lemma}{Lemma}
\newtheorem{problem}{Problem}
\newtheorem{ass}{Assumption}

\newcommand{\mc}{\mathcal}

\usepackage{enumitem}
\usepackage{xcolor} 
\usepackage{tikz} 
\usepackage{pgfplots} 
\usepgfplotslibrary{patchplots}
\usepackage{caption} 
\usepackage{subcaption} 
\usepackage{threeparttable} 

\usepackage{footnote}
\makesavenoteenv{tabular}
\makesavenoteenv{table}

\usepackage{comment}

\begin{document}
\begin{frontmatter}

\title{Data-Driven Behaviour Estimation in Parametric Games\thanksref{footnoteinfo}} 

\thanks[footnoteinfo]{The authors gratefully acknowledge financial support from ETHZ and the SNSF under the NCCR Automation grant 180545.}

\author[First]{Anna M. Maddux} 
\author[Second]{Nicol\`{o} Pagan} 
\author[Third]{Giuseppe Belgioioso}
\author[Third]{Florian Dörfler}

\address[First]{Sycamore, EPF Lausanne, Switzerland, (e-mail: anna.maddux@epfl.ch).}
\address[Second]{Social Computing Group, UZH, Switerzerland, (e-mail: nicolo.pagan@uzh.ch)}
\address[Third]{Automatic Control Laboratory, ETH Zürich, Switzerland, (e-mail: \{gbelgioioso,doerfler\}@ethz.ch)}

\begin{abstract}                
A central question in multi-agent strategic games deals with learning the underlying utilities driving the agents' behaviour.
%
Motivated by the increasing availability of large data-sets, we develop an unifying data-driven technique to estimate agents' utility functions from their observed behaviour, irrespective of whether the observations correspond to equilibrium configurations or to temporal sequences of action profiles.
%
%
Under standard assumptions on the parametrization of the utilities, the proposed inference method is computationally efficient and finds all the parameters that rationalize the observed behaviour best.
%
We numerically validate our theoretical findings on the market share estimation problem under advertising competition, using historical data from the Coca-Cola Company and Pepsi Inc. duopoly.
\end{abstract}


\end{frontmatter}

\section{Introduction}
The study of multi-agent systems (MAS) has been of interest to various research communities, such as economics, e.g. Bertrand--Nash competition \citep{narahari09}, computational social science, e.g. social networks formation
\citep{pagan19}, and automatic control, e.g. power systems \citep{belgioioso2022operationally}.
The decision-making processes in such MAS can be modelled as a strategic game, where, typically, the objective of each agent is to select an action that maximizes its own utility.

A pivotal problem in game theory deals with predicting the behaviour of agents from their utilities, namely, characterizing the game equilibria.
In many relevant cases, though, the underlying utilities are not directly observable. With the emergence of the internet of things, however, increasing amounts of data regarding the behaviour of each agent are available and can be used to infer the underlying utilities that rationally led to that behaviour. 
This line of research, termed \textit{inverse game theory}, has many useful applications: Inferred utilities can be used to predict the agents' behaviour in other similarly structured games, where the behaviour is unobservable \citep{waugh11}, or to cross-validate results from classical game theory, namely, to confirm whether a specified game admits any utilities that are consistent with the observed behaviour \citep{kuleshov15}.



Existing data-driven methods for inverse game theory can be broadly divided into two sub-classes: methods based on observing time-independent snapshots of a MAS, and methods based on action trajectories observations.
Among works of the former class, \citep{pagan19} and \citep{kuleshov15} propose an inverse-optimization approach to determine the utilities from an observed Nash and correlated equilibrium, respectively.
\citep{bertsimas15} combine ideas from inverse optimization and variational inequality (VI) theory to build a (non-) parametric estimation method.  
\citep{li19} and \citep{trivedi20} propose a data-driven framework specific to matchings and network formation, respectively.

The latter sub-class of methods (inverse-game-theoretic problems with action trajectories observations) has mostly been studied by the automatic control community. \citep{tsai16} and \citep{molloy17} propose an inference method based on the minimum principle of optimal control for estimating the unknown parameters of agents' cost functions from observed open-loop Nash equilibrium in two-player zero-sum and N-player dynamic games, respectively. 
\citep{inga19} apply the principle of maximum entropy in $N$-player inverse dynamic games
to infer the cost functions from both observed open-loop and feedback Nash equilibria.
Alternatively, \citep{nekipelov} develop a regret minimization method for auction games to infer players' private valuations (parameters). 

From the aforementioned literature, several critical issues emerge.
First, the objective is to find a parameter tuple which rationalizes the observed behaviour best\footnote{Under certain circumstances \citep{bertsimas15,pagan19} find sets (not singletons) of parameters.}.
However, the inverse game-theoretic problem is, in general, ill-posed, in the sense that possibly multiple parameter tuples exist which simultaneously are the most rational estimates of the observed behaviour. 
%
%
Moreover, all previously presented works are based on the assumption that the available observations are either snapshots or trajectories of action profiles, but always equilibrium configurations 
with the exception of \citep{nekipelov}.
In reality, however, interactions between agents may evolve over time, sometimes generating persistent fluctuations and equilibrium behaviour may not be observed. 
%

Motivated by these issues, in this paper, we develop a novel unifying data-driven behavior estimation method to learn parametric utilities in multi-agent games. The contributions are summarized next.
\begin{enumerate}
\item[(i)] Our inference method is unifying, in the sense that it handles observations of both snapshots and of trajectories of action profiles in a multi-agent game.
It builds on the assumptions that any observed snapshot corresponds to an approximate Nash equilibrium and that any observed action trajectory is compatible with an approximate better-response dynamics, i.e., in which each agent aims at improving its utility with respect to its previous action.
%
%

\item[(ii)]
Our method leverages a computationally tractable linear programming reformulation of the original inference problem to produce polyhedral sets containing all the parameters that rationalize the observed game behaviour best.
%
%
%
%
%
\item[(iii)]
We numerically validate our method 
on the market share evolution estimation problem under advertising competition \citep{gasmi91}, using historical data of advertising expenditures trajectories from the Coca-Cola Company and Pepsi Inc. duopoly. 
\end{enumerate}

\section{Mathematical setup}\label{sec:problem_formulation}
We consider a set of players (or agents) $\mathcal{I} :=\{1,\ldots,N\}$, where each agent $i \in \mathcal{I}$ chooses its action (or strategy) $a_i$ from its action set $\mc A_i $. Let $a:=(a_1,\ldots,a_N) \in \mc A$ be the action profile, where $\mc A:= \bigtimes_{i \in \mc I} \mc A_{i}$ is the global action set. Each agent $i \in \mathcal{I}$ aims at maximizing its own utility $U_i(a_i,a_{-i})$, which depends on its action $a_i\in\mathcal{A}_i$ and on the actions of the other agents $a_{-i}:=(a_1,\ldots,a_{i-1},a_{i+1},\ldots,a_N)\in\mathcal{A}_{-i}$, where $\mathcal{A}_{-i}:= \bigtimes_{j \in \mc I \setminus \{i\}  } \mc A_{j}$. The collection of players, strategy sets, and utility functions, i.e.,  ${\mathcal{G}:=(\mathcal{I},(\mathcal{A}_i)_{i\in\mathcal{I}},(U_i)_{i\in\mathcal{I}})}$, constitutes a game.

We assume that agents are rational decision makers. Namely, if an agent $i$ is allowed to revise its current action $\bar{a}_i$, it will choose an action $a_i \in \mathcal{A}_i$ that leads to a larger payoff, i.e., a so-called \textit{better response}.
\smallskip
\begin{definition}\label{def_better-response}
For all $i \in \mc I$, an action $a_i\in\mathcal{A}_i$ is an \textit{$\epsilon$-better response} to $\bar{a} = (\bar{a}_i, \bar{a}_{-i})\in\mathcal{A}_i\times\mathcal{A}_{-i}$, with $\epsilon\geq 0$, if
\begin{align}\label{better-response}
    U_i(a_i,\bar{a}_{-i})\geq U_i(\bar{a}_i,\bar{a}_{-i})-\epsilon.
\end{align}
Accordingly, the $\epsilon$-better response mapping is defined as
\begin{align*}
    \mathcal{R}_i^\epsilon(\bar{a})
    :=\{a_i\in\mathcal{A}_i\,~ \vert~ \eqref{def_better-response}\text{ holds}\}.
\end{align*}
If \eqref{better-response} holds with $\epsilon=0$, $a_i$ is known as better response.
{\hfill$\square$}
\end{definition}
An action $a_i$ is a \textit{best response} to $\bar{a}_{-i}$, if it is a better response compared to any alternative admissible action $\bar{a}_{i}\in\mathcal{A}_i$.
\smallskip
\begin{definition}\label{def_best-response}
For all $i \in \mc I$, an action $a_i\in\mathcal{A}_i$ is an \textit{$\epsilon$-best response} to $\bar{a}_{-i} \in \mc A_{-i}$, with $\epsilon \geq 0$, if
\begin{align}\label{best-response}
    U_i(a_i,\bar{a}_{-i})\geq U_i(\bar{a}_i,\bar{a}_{-i})-\epsilon, \quad \forall\bar{a}_i\in\mathcal{A}_i.
\end{align}
Accordingly, the $\epsilon$-best response mapping is defined as
\begin{align*}
    \mathcal{B}_i^\epsilon(\bar{a}_{-i}):=\{a_i\in\mathcal{A}_i\, ~\vert ~ \eqref{best-response}\text{ holds}\}.
\end{align*}
If \eqref{def_best-response} holds with $\epsilon=0$, $a_i$ is knwon as \textit{best response}.
{\hfill$\square$}
\end{definition}
\smallskip


If all agents play an $\epsilon$-best response action simultaneously, this action profile is an $\epsilon$-Nash equilibrium.

\smallskip
\begin{definition}\label{def_NE}
An action profile $ a^*=(a_1^*,\ldots,a_N^*)$ is an $\epsilon$-Nash equilibrium ($\epsilon$-NE), with $\epsilon \geq 0$, if
\begin{align}\label{NE}
\forall i\in\mathcal{I}:\quad 
a_i^* \in \mc B^\epsilon_i(a^*_{-i}).
\end{align}
If \eqref{NE} holds with $\epsilon=0$, $a^*$ is a Nash equilibrium (NE).
{\hfill$\square$}
\end{definition}

A NE is strategically-stable in the sense that no agent has an interest in unilaterally deviating. Next, we assume that the local utility functions admit a finite parametrization.

\begin{ass}\label{ass_parametric_utility}
For player $i \in \mathcal{I}$, the utility function $U_i(a_i,a_{-i};\theta_i)$ is parameterized by $\theta_i\in\Theta_i$, where $\Theta_i \neq \varnothing$ is the parameter space. Additionally, Definitions \ref{def_better-response}-\ref{def_NE} hold with respect to a specified parameterization.{\hfill$\square$}
\end{ass}


In this paper, we assume that the utilities are unknown but observations of the agents' actions are available. Thus, we focus on the problem of inferring the underlying parameters $\theta_i$ that led the system to the observed actions.
Concretely, we consider two types of observations:
\begin{enumerate}[label=(\roman*)]
    \item Dynamic observations: a temporal sequence of $m \in \mathbb{N}$ action profiles, ${(a(j))}_{j\in[m]}$, where $a(j)\in \mc A$ is the action profile at time step $j$, and $[m]:= \{1,\ldots,m \}$.

    \item Static observations: an equilibrium profile $a^* \in \mc A$. 
\end{enumerate}

We illustrate the difference with an example: in economics, competing firms repeatedly change their prices 
in response to other firms.
In this case, a sequence of action profiles, i.e., prices, is observed, and the firms' utility functions can be inferred by assuming that pairwise consecutive elements of the sequence are approximately compatible with better responses, i.e., $a_i(j\!+\!1)\in\mathcal{R}_i^{\epsilon_i}(a(j))$, for some $\epsilon_i \geq 0$.
On the other hand, if the firms have reached an equilibrium, their utility functions can still be inferred by assuming that the observed action profile $a^*$ is approximately a NE, namely, $a_i^* \in \mc B_i^{\epsilon_i}(a_{-i}^*)$ for all $i \in \mc I$.
The term ``approximately'' stems from the fact that the agents may not act fully rational. Namely, the observed action trajectories and equilibrium profiles may violate some of the perfect better-response \eqref{better-response} and NE conditions \eqref{def_NE}, respectively.

\section{Data-driven behavior estimation}
Before proceeding, this next assumption helps to unify the notation for the static and dynamic cases.
\begin{ass}
\label{ass:BRNEcond}
For all $i\in\mathcal{I}$, we construct a set of $m$ data points $\mc D_i:=\{(a_i^j,\bar{a}^j)\}_{j \in [m]}$, where $a_i^j \in \mc A_i$ and $\bar{a}^j =(\bar{a}_i^j,\bar{a}_{-i}^j) \in \mc A$.
\begin{enumerate}[label=(\roman*)]
    \item If the observations are dynamic, the data points are of the form
$(a_i^j,\bar{a}^j)=(a_i(j\!+\!1),a(j))$, for all time steps $j\in\left[m\right]$, and $a_i(j+1)\in\mathcal{R}_i^{\epsilon_i}(a(j))$, for all $j\in\left[m\right]$, w.r.t. the true parameterization $\theta_i$ of $U_i(\cdot,\cdot;\theta_i)$ for some $\epsilon_i \geq 0$. 
    \item If the observations are static, the data points are of the form $(a_i^j,\bar{a}^j)\!=\!\left(a_i^*,(\tilde a_i^j, a_{-i}^*)\right)$, where $\tilde a_i^j \!\in\! \mc A_i$, for all $j\!\in\!\left[m\right]$, are $m$ alternative admissible actions of player $i$\footnote{The equilibrium profile $(a_i^*,a_{-i}^*)$ is observed while $\tilde{a}_i^j\in\mathcal{A}_i$, $j\in\left[m\right]$, are fictitious admissible actions agent $i$ could have played. 
    }.
    Furthermore, $a_i^* \in \mc B_i^{\epsilon_i}(a_{-i}^*)$ w.r.t. to the true parameterization $\theta_i$ of $U_i(\cdot,\cdot;\theta_i)$ for some $\epsilon_i \geq 0$.
    
{\hfill $\square$}
\end{enumerate} 
\end{ass}
This assumption states that any observed action trajectory $(a_i(j))_{j\in[m]}$ is an $\epsilon$-better response dynamics and any observed action profile $a^*$ is an $\epsilon$-NE for the true paramerization of the utility function.
Then, our objective is to find all admissible preference vectors $\theta_i$ that fully rationalize the observed behaviour, namely, the action trajectories or the NE. Alternatively, if no such vector exists, we aim to find all the vectors $\theta_i$ that violate the better response or NE conditions the least.
To measure the irrationality, we introduce the error functions $e_i: \mc A_i \times \mc A \times \Theta_i \rightarrow \mathbb R$, defined as in \citep{pagan19}
\begin{align}\label{pre_error_func}
    e_i(a_i,\bar{a};\theta_i):=U_i(\bar{a}_i,\bar{a}_{-i};\theta_i)-U_i(a_i,\bar{a}_{-i};\theta_i),
\end{align}
which takes on positive values whenever agent $i$'s action $a_i$ yields a lower utility than $\bar{a}_i$, given that the other agents' actions $\bar{a}_{-i}$ remains fixed. 
For dynamic observations, $\theta_i$ is a rational preference if and only if $a_i(j\!+\!1)$ leads to higher utility compared to $a_i(j)$, when $a_{-i}(j)$ remains fixed. 
For static observations, instead, we expect that $a_i^*$ leads to higher payoffs compared to any other admissible alternative $\tilde a_i \in \mathcal{A}_i$. 
Regardless of whether the observations are evolving or static, 
if  $e_i(a_i,\bar{a};\theta_i)\leq 0$ then $\theta_i$ is a rational behavior. Vice-versa,
the irrationality of $\theta_i$ comes from the positive contributions of the error function in \eqref{pre_error_func}: $e_i^+(a_i,\bar{a};\theta_i):=\max\{0,e_i(a_i,\bar{a};\theta_i)\}$.

To measure an agent's maximal violation of the rationality assumption, 
we introduce a function,
\begin{align}\label{distance_func}
    d_i\left(\theta_i; \mc D_i \right):=\max_{j\in\left[m\right]}\, e_i^+(a_i^j,\bar{a}^j;\theta_i),
\end{align}
termed  the irrationality loss function. Any preference vector $\hat{\theta}_i$ which minimizes $ d_i\left(\, \cdot \,; \mc D_i \right)$ ensures that the worst violation of the rationality assumption over all $m$ observations $\mc D_i$ is minimal and $\hat{\theta}_i$ rationalizes the observed individual behaviour best. 
Hence, our inference problem can be cast as $N$ decoupled optimization problems:
\smallskip\noindent
\begin{problem}\label{inference_problem}
Let assumptions \ref{ass_parametric_utility}-\ref{ass:BRNEcond} hold. For all $i \in \mc I$, find all the vectors of preferences $\hat{\theta}_i$ that solve
\begin{align} \label{eq:pr1}
    \left\{
    \begin{array}{r l}
    \displaystyle
\min_{\theta_i  \in \mathbb{R}^{p_i}} &  d_i(\theta_i ; \mc D_i)\\
\text{s.t.} & \theta_i \in \Theta_i.
\end{array}    
\right.
\end{align}

\end{problem}


Next, we show that the observed trajectory of action profiles $(a(j))_{j\in\left[m\right]}$ can be interpreted as an $\epsilon$-better response dynamics for each agent $i$ and, if the available data points $\mc D_i$ in Assumption~\ref{ass:BRNEcond}~(ii) are exhaustive, i.e., $\mc D_i$ contains all feasible actions of agent $i$, the observed static profile $a^*$ can be interpreted as an $\epsilon$-NE\footnote{Note that we make no assumptions about the utility functions w.r.t. the actions. Thus, if the game is ill-posed, i.e., a NE does not exist, then our method finds a preference vector that violates the NE conditions the least.}, where $\epsilon$ is the maximum value among the solutions of the $N$ problems in \eqref{eq:pr1}.

\begin{prop} Let assumptions \ref{ass_parametric_utility}-\ref{ass:BRNEcond} hold true. \label{prop:1}
\begin{enumerate}[label=(\roman*)]
\item For $i\in \mathcal{I}$, the observed action trajectory, $(a_i(j))_{j\in\left[m\right]}$ is an $\bar{\epsilon}$-better response dynamics w.r.t. $U_i(\cdot,\cdot,\hat{\theta}_i)$, where $\bar{\epsilon}:=\max_{i\in\mathcal{I}} d_i(\hat{\theta}_i; \mc D_i)$, and $\hat \theta_i$ as in \eqref{eq:pr1}.
\item Assume that $\bigcup_{i=1}^m \tilde{a}_i^j=\mc A_i$, for all $i \in \mathcal{I}$. Then, the observed action profile $a^*$ is an $\bar \epsilon$-NE w.r.t. $U_i(\cdot,\cdot,\hat{\theta}_i)$, where $\bar \epsilon:=\max_{i\in\mathcal{I}} d_i(\hat{\theta}_i; \mc D_i)$, and $\hat \theta_i$ as in \eqref{eq:pr1} for all $i \in \mathcal{I}$.
 {\hfill$\square$}\\
\end{enumerate}
\end{prop}
For the proof we refer the reader to Appendix \ref{ap:proofOfProp1}.

%

\subsection{Behavior estimation in contextual games}
In this section, we extend the formulation of Problem~\ref{inference_problem} to the case of multiple data-sets, i.e., multiple trajectories or equilibria observations.
In many relevant scenarios, (the outcome of) multi-agent games depend on exogenous but observable factors (e.g. GDP, shocks, seasonality in Bertrand--Nash games).
%
More formally, we denote such factors as the context $\xi\in\Xi$ of the game, where $\Xi\subseteq \mathbb R^c$ is the set of possible contexts \citep{sessa20}. Hereafter, we assume that each agent's parametric utility depends also on a known context $\xi$, as formalized next.
\begin{ass}
For all $i\in\mathcal{I}$, the local utility function $U_i(a_i,a_{-i};\theta_i,\xi)$ depends also on an observable parameter (context) $\xi\in \Xi$.
{\hfill $\square$}
\end{ass}


To generalize the inference problem in Section~\ref{sec:problem_formulation}, we assume to have $n$  observations, i.e., action profile trajectories $(a^k(j))_{j \in [m^k]}$ or equilibrium profiles $a^{*,k}$, each connected to a different game context $\xi^k$, with $k \in [n]$. As in Section~\ref{sec:problem_formulation}, we use the following assumption to construct a data-set that helps to unify the notation for the static and dynamic cases.
\begin{ass}
\label{ass:multi-dataSet}
Let $\{ \xi^k \}_{k \in [n]}$, with $\xi^k \in \Xi$, be a set of contexts of a multi-agent game ${\mathcal{G}=(\mathcal{I},(\mathcal{A}_i)_{i\in\mathcal{I}},(U_i)_{i\in\mathcal{I}},\Xi)}$. For each context $\xi^k$, we construct $n$ sets of data points $\{ \mc D_i^k\}_{i \in \mathcal{I}}$ as in Assumption~\ref{ass:BRNEcond}, where $D_i^k:=${\small $\{(a_i^{j,k},\bar{a}^{j,k})\}_{j \in [m^k]}$} and $m^k$ is the number of data points relative to the $k$-th context.
{\hfill $\square$}
\end{ass}
The multi data-sets extension of the original behaviour  estimation problem (Problem \ref{inference_problem}) can be recast as $N$ decoupled optimization problems, as formalized next.

\begin{problem}\label{context_inference_problem}
Let assumptions~\ref{ass_parametric_utility}-\ref{ass:multi-dataSet} hold true. For all $i \in \mathcal{I}$, find all vectors of preferences $\hat{\theta}_i$ that solve
\begin{align} \label{eq:pr2}
   \left\{
    \begin{array}{r l}
    \displaystyle
\min_{\theta_i \in \mathbb{R}^{p_i}} & \displaystyle \max_{k \in [n]}\;  d_i(\theta_i ; \mc D^k_i, \xi^k)\\
\text{s.t.} & \theta_i \in \Theta_i,
\end{array}    
\right.
\end{align}
where the distance function\footnote{With some abuse of notation, we keep the same notation for the functions $e_i$, $e_i^+$ and $d_i$ in contextual games%
.} $d_i$ in \eqref{eq:pr2} is now parametric also with respect to the game context $\xi^k$, i.e.,
\begin{align}\label{context_distance_func}
    d_i(\theta; \mc D^k_i, \xi^k):=\max_{j\in\left[m^k\right]}e_i^+(a_i^{j,k},\bar{a}^{j,k};\theta_i,\xi^k).
\end{align}

\end{problem}


\subsection{A linear programming approach}\label{sec:LP_solution}
Before proceeding, we assume that the utilities are linearly parametrized, as in \citep{waugh11,kuleshov15,pagan19}.

\begin{ass}
\label{eq:linParam}
For all $i \in \mathcal{I}$, the function $U_i(a_i,a_{-i};\theta_i,\xi)$ is linear in its parameters $\theta_i\in\Theta_i$.{\hfill$\square$}
\end{ass}

It follows that $d_i$ in \eqref{context_distance_func} is piece-wise linear, thus convex, in its parameters $\theta_i$. Therefore, the minimax optimization in Problem~\ref{context_inference_problem} can be equivalently recast as a linear program (LP), see e.g.  \S~1 in \citep{ahuja85}, which is always feasible.
\smallskip\noindent
\begin{prop}\label{LP_reformulation}
Let the preconditions of Problem~\ref{context_inference_problem} hold. For all $i\! \in\! \mathcal{I}$, vector $\hat \theta_i$ solves \eqref{eq:pr2} if and only if there exists a constant $\hat \epsilon_i$ such that the pair $(\hat{\theta}_i,\hat \epsilon_i)$ solves the LP
\begin{equation}
\label{eq:pr3}
\left\{
\begin{array}{r l}
\displaystyle
\min_{\epsilon_i ,\theta_i }
& \epsilon_i\\[.5em]
\text{s.t.} & e_i(a_i^{j,k},\bar{a}^{j,k};\theta_i,\xi^k)\leq\epsilon_i, \; \forall  j\!\in \![m^k],\, \forall k \in [n],\\[.5em]
      &\epsilon_i\geq 0,\\[.5em]
      &\theta_i\in\Theta_i.  
\end{array}
\right.
\end{equation}
%

\vspace*{-.7em}
{\hfill $\square$}
\end{prop}
%
\begin{lemma}
For all $i \in \mathcal{I}$, the solution set of the LP in \eqref{eq:pr3} is nonempty. {\hfill$\square$}
\end{lemma}
For the proof we refer the reader to Appendix \ref{ap:ProofOfLPref}.

Finally, by leveraging the LP reformulation in Proposition~\ref{LP_reformulation}, we can characterize the solution set of the original multi data-set behavior estimation problem (Problem~\ref{context_inference_problem}).

\begin{prop}\label{prop:sol}
Let assumptions \ref{ass_parametric_utility}-\ref{ass:multi-dataSet} hold true. For all $i \in \mc I$, the solution set $\hat{\Theta}_i$ of \eqref{eq:pr2} reads as
\begin{align}\label{solution_set}
\begin{split}
\hat{\Theta}_i:=\{&\theta_i\in\Theta_i\, ~|~ e_i(a_i^{j,k},\bar{a}^{j,k};\xi^k,\theta_i)\leq \hat{\epsilon}_i,\\
&\forall j\in [m^k],\ \forall k\in\left[n\right]\Big\},
\end{split}
\end{align}
where $\hat{\epsilon}_i$ is the optimal value of the LP in \eqref{eq:pr3}.
{\hfill $\square$}
\end{prop}
Note that $\hat{\Theta}_i$ is a polyhedron which contains all parameters $\theta_i$ that rationalize best the observed behaviour of agent $i$ over all contexts $\xi^k$. If $\hat\epsilon_i\!=\!0$, then $\max_{k \in [n]}\;  d_i(\hat \theta_i ; \mc D^k_i, \xi^k)\!=\!0$ for all $\hat \theta_i\in\hat{\Theta}_i$, which implies that agent $i$'s observed behaviour is perfectly rationalizable.

\section{Numerical and empirical analysis}\label{sec:numerical_analysis}
In this section, we numerically validate our inference method on (a) the demand estimation problem under Bertrand--Nash competition, using static observations of NE prices, and on (b) the market share evolution estimation problem under advertising competition \citep{gasmi91}, using dynamic observations constituted by historical data of advertising expenditures.

\subsection{Demand estimation in Bertrand--Nash competition}\label{static_application}

Consider two firms competing to maximize their revenues by setting prices $p_1$ and $p_2$, respectively. The demand of each firm $i\in\{1,2\}$ is modelled by the parametric function
\begin{align}\label{demand_func}
D_i(p_1,p_2;\theta_i,\xi)=\theta_{i,0}+\theta_{i,1}p_1+\theta_{i,2}p_2+\theta_{i,3}\xi,
\end{align}
where $\xi$ represents time and region specific economic indicators (e.g. shocks, seasonality). We assume $\theta_{i,i}\leq 0$ and $\theta_{i,-i}\geq 0$ such that, for firm $i$, an increase of price $p_i$ has a negative effect, and an increase of price $p_{-i}$ has a positive effect on its demand ($\theta_{i,-i}=\theta_{1,2}$, if $i=1$ and, $\theta_{i,-i}=\theta_{2,1}$, if $i=2$). Then, given context $\xi$, each firm sets its prices to maximize its revenue, i.e., $U_i(p_1,p_2;\theta_i,\xi)=p_i\cdot D_i(p_1,p_2;\theta_i,\xi)$.

Given a set of observed contexts $\{\xi^k\}_{k\in\left[n\right]}$, we seek to estimate the most rational normalized\footnote{Normalization ensures that the everywhere-zero revenue function does not let all prices appear rational \citep{ng00}.} parameters {\small$\Tilde{\theta}_i:=(\frac{\theta_{i,0}}{\theta_{i,-i}}, \frac{\theta_{i,i}}{\theta_{i,-i}},\frac{\theta_{i,3}}{\theta_{i,-i}})^\top$} from a set of observed price tuples $\{(p_1^{*,k},p_2^{*,k})\}_{k\in\left[n\right]}$, each connected to a different context $\xi^k$ and assumed to be a NE. Thus, we construct $n$ batches of data-sets $D_i^k:=\{(p_i^{*,k},\tilde{p}_i^{j,k},p_{-i}^{*,k})\}_{j \in [m^k]}$, where $\{(\tilde{p}_1^{j,k},\tilde{p}_2^{j,k})\}_{j\in m^k}$ is a set of alternative feasible price profiles. 
This demand estimation problem
satisfies the assumptions of our framework. Thus, its solutions can be found by exploiting Proposition~\ref{prop:sol}.

For our numerical study, we assume the true parameters $\theta_1^*=(1, -1.2, 0.5, 1)^\top$ and $\theta_2^*=(1, 0.3, -1, 1)^\top$ of the demand function \eqref{demand_func} to be unknown. Drawing inspiration from \S~8.1 in \citep{bertsimas15}, we model contexts $\xi^k\in\{\xi^k\}_{k\in\left[50\right]}$ as i.i.d. normally distributed random variables with mean 5 and standard deviation 1.5. and solve for the NE prices $(p_1^{*,k},p_2^{*,k})$ connected to context $\xi^k$. For each agent  $i$ and each context $\xi^k$, we estimate $\Tilde{\theta}_i$ from data-set\footnote{$\tilde{p}_i^{j,k}$ is the $j$-th point in a regular  price grid $\left[0,p_{\max}\right]$ of $2^7+1$ points, which can be seen as a discretization of the price sets $\mathcal{A}_i$} $D_i^k:=\{(p_i^{*,k}, \tilde{p}_i^{j,k},p_{-i}^{*,k})\}_{j\in\left[2^7+1\right]}$.
To further validate our method, we compare its results with those obtained via the VI-based method proposed in \S~3.1 of \citep{bertsimas15}. 
The polyhedral solution set $\hat{\Theta}_1$ of firm 1 contains $\theta_1^*$ and $\theta_1^{\text{VI}}$. 
Since for each context $\xi^k$, the price pair $(p_1^{*,k},p_2^{*,k})$ is a NE, both methods unsurprisingly find parameter vectors that perfectly rationalize the observed prices. However, the solution of the inverse VI problem $\theta_1^{\text{VI}}$ does not coincide with the true parameter vector $\theta_1^*$. On the other hand, solving Problem \ref{context_inference_problem} gives all the parameter vectors that perfectly rationalize the generated prices, including the true parameter vector $\theta_1^*$.

\begin{figure}[t]
         \centering
         \input{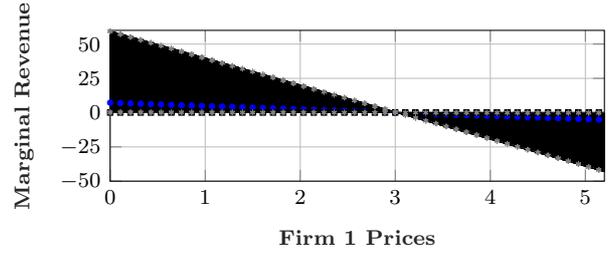}
         \caption{The true marginal revenue function (blue line), the inverse-VI fit \citep{bertsimas15} (black line, square markers), and the fit obtained from the polyhedral solution set which reconcile all the data exactly (grey region) for both firms.}
         \label{fig:marginal_revenue}
     
\end{figure}
We plot the fitted marginal revenue function $\text{MR}_1=\frac{\partial}{\partial p_i}p_i\cdot D_1(p_1,p_2;\theta_i,\xi)$ of firm 1 obtained from the vertices of the polyhedral solution sets $\hat{\Theta}_1$, see Figure~\ref{fig:marginal_revenue}. The grey area displays all potential marginal revenue functions in the parametric family specified in \eqref{demand_func} which exactly reconcile the simulated NE prices. Note that it contains the true marginal revenue function and the fitted function obtained from \citep{bertsimas15} since both reconcile the data. 
%


\begin{table}
\captionsetup{font=small}
\caption{
Rationalization error $\epsilon_i$ and estimated beliefs $k_{1,i},k_{2,i}$ obtained via our method ($L_\infty$), the $L_2$-based variant, and OLS regression.}
\label{tab:L_infty_L2_ols}
\centering
\begin{tabular}{|c||c|c|c|c|c|c|}
\hline
& \multicolumn{3}{|c|}{Coke} & \multicolumn{3}{|c|}{Pepsi}\\
\hline
& $\hat{\epsilon}_1$ & $\hat{k}_{\{1,1\}}$ & $\hat{k}_{\{2,1\}}$ &  $\hat{\epsilon}_2$ & $\hat{k}_{\{1,2\}}$ & $\hat{k}_{\{2,2\}}$ \\
\hline
$L_\infty$ & 0.018 & 0.066 & 0.057 & 0.011 & 0.036 & 0.022 \\
\hline
 $L_2$ & 0.021 & 0.078 & 0.052 & 0.013 & 0.021 & 0.016 \\
\hline
 OLS & 0.059 & 0.011 & 0.007 & 0.015 & 0.012 & 0.011\\
\hline
\end{tabular}
\end{table}

\subsection{Estimating market shares in advertising competition}\label{dynamic_application}
In this subsection, we formulate the problem of estimating the beliefs of two firms in an advertising competition on the market share evolution. Under the assumption that the trajectory of observed advertising expenditure profiles $(a_1(j),a_2(t))_{t\in[T]}$ is compatible with an approximate better-response dynamics (Assumption \ref{ass:BRNEcond}), we can solve the belief estimation problem using Proposition~\ref{prop:sol}.

We model the market share evolution via the parametric Lanchester model (LM) of advertising competition in  \citep{chintagunta92}. Its discrete-time specification, taken from Eqn.~(2) in \citep{chintagunta92}, is given by
\begin{align}\label{lanchester}
\begin{split}
    M(t\!+\!1)=M(t)&+k_1\sqrt{a_1(t\!+\!1)}(1-M(t))\\
    &-k_2\sqrt{a_2(t\!+\!1)}M(t),
\end{split}
\end{align}
where $M$ is the market share of firm~1, $(1\!-\!M)$ is the share of firm~2, and $k_1, k_2\geq 0$ are measures of the advertising effectiveness of firm~1 and firm~2, respectively.
To include word-of-mouth effects on the market share evolution, \citep{sorger89}  suggests a modification of the LM, as shown in Eqn.~(4) of \citep{chintagunta95}
\begin{align}\label{sorger}
\begin{split}
    M(t\!+\!1)=M(t)&+k_1\sqrt{a_1(t\!+\!1)}\sqrt{(1-M(t))}\\
    &-k_2\sqrt{a_2(t\!+\!1)}\sqrt{M(t)}.
\end{split}
\end{align}

In a duopolistic market share competition, both firms update their advertising expenditure to gain market share. Moreover, each firms's choice of advertising expenditure depends on its beliefs on its own and the competitor's advertising effectiveness.
Let us denote by $\theta_1=(k_{1,1},k_{2,1})$ and $\theta_2=(k_{1,2},k_{2,2})$ firm~1's and firm~2's beliefs on their true advertising effectiveness $k_1$ and $k_2$, respectively.
We assume that firm~1 sets its new advertising expenditure, $a_1(\!t+1\!)$, with the objective of increasing its current market share, i.e., $M(a_1(\!t+1\!),a_2(\!t\!);\theta_1) \geq M(a_1(t),a_2(t);\theta_1)$. Similarly, we assume that firm~2 sets $a_2(\!t+1\!)$ such that $1-M(a_1(t),a_2(t\!+\!1);\theta_2) \geq 1 -M(a_1(t),a_2(t);\theta_2)$. Thus, $\mc D_i:=\{(a_i(\!t+\!1),a_i(t),a_{-i}(t))\}_{t\in [T]}$, constitutes the data-set of agent $i$, from which we estimate $\hat{\theta}_i$ such that the observed advertising expenditure trajectory is compatible with a better-response dynamics. 


\subsubsection{Empirical analysis} \label{sec:empircal_analysis}
In this section, we show the applicability of our method to a real-world data-set.
The Coca-Cola Company (firm 1) and Pepsi Inc. (firm 2) are amongst the most famous examples of duopolistic rivalries over market share, while other soft-drinks companies own a negligible fraction of market shares. We assume Coke and Pepsi set their advertising expenditures to gain market share, where Coke and Pepsi consider the market share to evolve according to the LM \eqref{lanchester} and SM, respectively. Thus, using data\footnote{The description of the data-set, generously provided by Gasmi, Laffont and Vuong, can be found in \citep{gasmi91}.} from 1968 to 1986, which spans a period of intense advertising competition, we estimate Coke's and Pepsi's beliefs on the market share evolution or, more concretely, the advertising effectiveness parameters using Proposition \ref{prop:sol}.

We find\footnote{The code is made available at https://github.com/amaddux9/\\Data-Driven-Behaviour-Estimation.} polyhedral parameter sets of Coke and Pepsi which rationalize the observed advertising expenditures up to an error of at most $\hat{\epsilon}_i$, $i\in\{1,2\}$, reported in Table~\ref{tab:L_infty_L2_ols}. Coke's and Pepsi's estimated effectiveness parameters $\hat{\theta}_1=(\hat{k}_{1,1},\hat{k}_{2,1})$ and $\hat{\theta}_2=(\hat{k}_{1,2},\hat{k}_{2,2})$, respectively, are positive, as required, and the believed  advertising effectiveness of Coke is slightly larger than that of Pepsi, which agrees with qualitative findings \citep{yoffie06}. 
%
%
In Figure~\ref{fig:irrationality} we compare our results with those obtained using a $L_2$-modification of our proposed method\footnote{The irrationality loss function is defined as $d_i(\theta_i,\mc{D}_i)=(1/m\sum_{j\in[m]}e:i^+(a_i^j,\overline{a}^j;\theta_i)^2)^{\frac{1}{2}}$} and the ordinary least square method (OLS), which differs fundamentally from the other two as it simply performs an ordinary least squares regression on the actual market share evolution.
As expected, our method minimizes the maximum deviation from rational better-responding behavior, while OLS regression leads to a much more irrational behavior (since it neglects any relationship between consecutive action pairs).
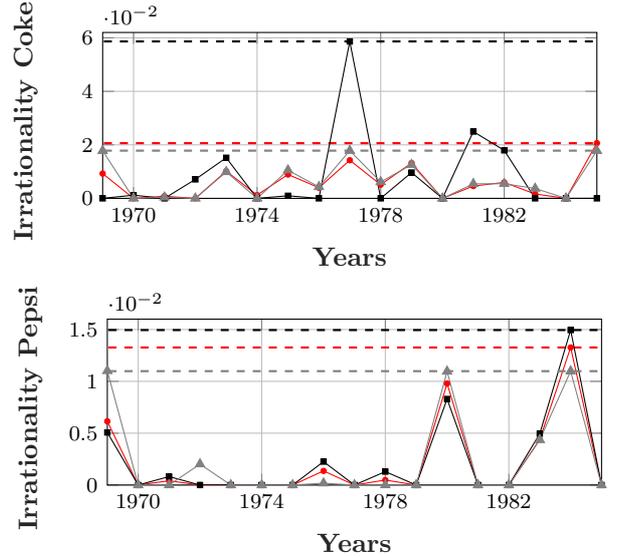
\begin{figure}[t!]
\captionsetup{font=footnotesize}
     \centering
     \begin{subfigure}{0.45\textwidth}
         \centering
%
%
\begin{tikzpicture}
\pgfplotsset{every tick label/.append style={font=\small}}

\begin{axis}[%
width=6.5cm,
height=2.2cm,
at={(0.758in,0.481in)},
scale only axis,
xmin=1,
xmax=17,
xtick={2,6,10,14},
xticklabels={{1970},{1974},{1978},{1982}},
xlabel style={font=\bfseries\color{white!15!black}},
xlabel={Years},
ymin=0,
ymax=0.062,
ylabel style={font=\bfseries\color{white!15!black},shift={(-0.075,-0.1)}},
ylabel={Irrationality Coke},
axis background/.style={fill=white},
xmajorgrids,
ymajorgrids
]

\addplot [color=red, mark size=1.0pt, mark=*, mark options={solid, fill=red, red}, forget plot]
  table[row sep=crcr]{%
1	0.00925766016669168\\
2	0\\
3	0.000791433419991438\\
4	0\\
5	0.0101703878243038\\
6	0.00122451838666127\\
7	0.00892697849087411\\
8	0.00398556662834304\\
9	0.0142443855939668\\
10	0.00511795179323508\\
11	0.0132013989143733\\
12	0\\
13	0.00459863637248018\\
14	0.00596264299558769\\
15	0.0016473384173516\\
16	0\\
17	0.02063998840364\\
};
\addplot [domain = 0:20,
        thick,
        dashed,
        red,
        ]{0.02063998840364};
\addplot [color=black, mark size=1.0pt, mark=square*, mark options={solid, fill=black, black}, forget plot]
  table[row sep=crcr]{%
1	0\\
2	0.00119562671570045\\
3	0\\
4	0.00708791767429975\\
5	0.0151735516737285\\
6	0\\
7	0.000967077740188358\\
8	0\\
9	0.0586465936879948\\
10	0\\
11	0.00966224498051532\\
12	0\\
13	0.0249567562643033\\
14	0.0179442067235416\\
15	0\\
16	0\\
17	0\\
};
\addplot [domain = 0:20,
        thick,
        dashed,
        black,
        ]{0.0586465936879948};
        
\addplot [color=gray, mark size=2.0pt, mark=triangle*, mark options={solid, fill=gray, gray}, forget plot]
  table[row sep=crcr]{%
1	0.0177850991488373\\
2	0\\
3	0.000402745007615294\\
4	0\\
5	0.0098844854426259\\
6	0\\
7	0.0105208091634327\\
8	0.00429177221663265\\
9	0.0178283670186901\\
10	0.00613553548598262\\
11	0.0125868404255104\\
12	0\\
13	0.0054252369597059\\
14	0.00547167137856165\\
15	0.00369313326748157\\
16	0\\
17	0.017816617125303\\
};
\addplot [domain = 0:20,
        thick,
        dashed,
        gray,
        ]{0.0178283670186901};
\end{axis}
\end{tikzpicture}
     \end{subfigure}
     \begin{subfigure}{0.45\textwidth}
         \centering
%
%
\begin{tikzpicture}
\pgfplotsset{every tick label/.append style={font=\small}}

\begin{axis}[%
width=6.5cm,
height=2.2cm,
at={(3.327in,0.481in)},
scale only axis,
xmin=1,
xmax=17,
xtick={2,6,10,14},
xticklabels={{1970},{1974},{1978},{1982}},
xlabel style={font=\bfseries\color{white!15!black}},
xlabel={Years},
ymin=0,
ymax=0.016,
ylabel style={font=\bfseries\color{white!15!black},shift={(-0.075,-0.1)}},
ylabel={Irrationality Pepsi},
axis background/.style={fill=white},
xmajorgrids,
ymajorgrids
]

\addplot [color=red, mark size=1.0pt, mark=*, mark options={solid, fill=red, red}, forget plot]
  table[row sep=crcr]{%
1	0.00614621792185613\\
2	0\\
3	0.000412756688014139\\
4	0\\
5	0\\
6	0\\
7	0\\
8	0.00136454287504726\\
9	0\\
10	0.000483060039724137\\
11	0\\
12	0.0097911042347349\\
13	0\\
14	0\\
15	0.0043942905074434\\
16	0.0132666940230387\\
17	0\\
};
\addplot [domain = 0:20,
        thick,
        dashed,
        red,
        ]{0.0132666940230387};
\addplot [color=black, mark size=1.0pt, mark=square*, mark options={solid, fill=black, black}, forget plot]
  table[row sep=crcr]{%
1	0.00506938990468067\\
2	0\\
3	0.000817249143086333\\
4	0\\
5	0\\
6	0\\
7	0\\
8	0.00225715804507887\\
9	0\\
10	0.00129442435388853\\
11	0\\
12	0.00830025296705543\\
13	0\\
14	0\\
15	0.00495367132561757\\
16	0.0149469420083142\\
17	0\\
};
\addplot [domain = 0:20,
        thick,
        dashed,
        black,
        ]{0.0149469420083142};

\addplot [color=gray, mark size=2.0pt, mark=triangle*, mark options={solid, fill=gray,gray}, forget plot]
  table[row sep=crcr]{%
1	0.0110104072580773\\
2	0\\
3	0\\
4	0.00202075111718122\\
5	0\\
6	0\\
7	0\\
8	0.000176488851323747\\
9	0\\
10	0\\
11	0\\
12	0.0109584201121063\\
13	0\\
14	0\\
15	0.00435469042103785\\
16	0.0109833271217873\\
17	0\\
};
\addplot [domain = 0:20,
        thick,
        dashed,
        gray,
        ]{0.0109833271217873};
\end{axis}
\end{tikzpicture}%
     \end{subfigure}
     \caption{Coke's and Pepsi's irrationality, expressed by violations of the better-response conditions in \eqref{better-response}, over the time horizon. Resulting irrationality from the OLS regression fit (black line), the $L_2$-based variant (red line), and from our estimation method (grey line).}
  \label{fig:irrationality}
\end{figure}
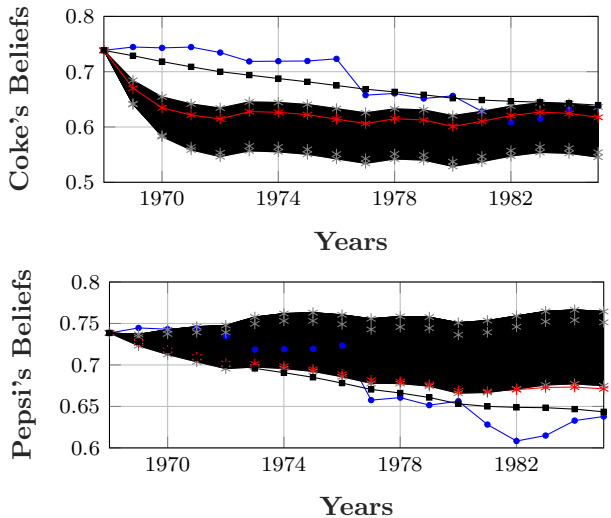
\begin{figure}[!htb]
\captionsetup{font=footnotesize}
    \centering
     \begin{subfigure}{0.45\textwidth}
         \centering
%
%
\begin{tikzpicture}
\pgfplotsset{every tick label/.append style={font=\small}}

\begin{axis}[%
width=6.5cm,
height=2.2cm,
at={(0.758in,0.481in)},
scale only axis,
xmin=1,
xmax=18,
xtick={3,7,11,15},
xticklabels={{1970},{1974},{1978},{1982}},
xlabel style={font=\bfseries\color{white!15!black}},
xlabel={Years},
ymin=0.5,
ymax=0.8,
ylabel style={font=\bfseries\color{white!15!black},shift={(0.0,-0.05)}},
ylabel={Coke's Beliefs},
axis background/.style={fill=white},
xmajorgrids,
ymajorgrids,
]
\addplot [color=blue, mark size=1.0pt, mark=*, mark options={solid, fill=blue, blue}, forget plot]
  table[row sep=crcr]{%
1	0.738591230987051\\
2	0.744620276849233\\
3	0.743087144080133\\
4	0.744480014058034\\
5	0.734562817374239\\
6	0.718587384502697\\
7	0.719016091196367\\
8	0.719341418235089\\
9	0.723257351397837\\
10	0.657583595317576\\
11	0.660627889831085\\
12	0.651520953831866\\
13	0.656345693598383\\
14	0.628155103863401\\
15	0.608293382967169\\
16	0.614894888824271\\
17	0.632923740112681\\
18	0.637850059562571\\
19	0.624847444924705\\
};
\addplot [color=gray, mark size=2.0pt, mark=asterisk, mark options={solid, gray}, forget plot]
  table[row sep=crcr]{%
1	0.738591230987051\\
2	0.684317304955061\\
3	0.654432223312138\\
4	0.641690498312399\\
5	0.634257867592339\\
6	0.645309120870074\\
7	0.644382127056871\\
8	0.640689923105775\\
9	0.633228873171658\\
10	0.625645518520773\\
11	0.632805708756496\\
12	0.630987943398583\\
13	0.620549873532067\\
14	0.628340427514512\\
15	0.63764977484573\\
16	0.643802316473428\\
17	0.642386541028871\\
18	0.63608742479616\\
19	0.612536647258844\\
};
\addplot [color=gray, mark size=2.0pt, mark=asterisk, mark options={solid, gray}, forget plot]
  table[row sep=crcr]{%
1	0.738591230987051\\
2	0.64100994859013\\
3	0.584566456075738\\
4	0.563567207208195\\
5	0.553250337440242\\
6	0.565560852285104\\
7	0.563872705858414\\
8	0.559432381277899\\
9	0.550955977434403\\
10	0.542675416804365\\
11	0.551401790355387\\
12	0.54915980808043\\
13	0.537262898611903\\
14	0.546669488720833\\
15	0.557150469853094\\
16	0.563731980589134\\
17	0.561606952528243\\
18	0.554240137577236\\
19	0.527501396432295\\
};
\addplot [color=gray, mark size=2.0pt, mark=asterisk, mark options={solid, gray}, forget plot]
  table[row sep=crcr]{%
1	0.738591230987051\\
2	0.642112454619908\\
3	0.582558830035727\\
4	0.558526896543968\\
5	0.546393596811143\\
6	0.556757482334975\\
7	0.555002456237832\\
8	0.550671614211642\\
9	0.542429379679178\\
10	0.534219092369346\\
11	0.542083183957094\\
12	0.540073660176655\\
13	0.528750629421901\\
14	0.537310322087951\\
15	0.547380014241206\\
16	0.553985910149654\\
17	0.55234943282698\\
18	0.545452425923719\\
19	0.519883315557908\\
};

\addplot[area legend, draw=black, fill=black, fill opacity=0.2, forget plot]
table[row sep=crcr] {%
x	y\\
1	0.738591230987051\\
2	0.64100994859013\\
3	0.584566456075738\\
4	0.563567207208195\\
5	0.553250337440242\\
6	0.565560852285104\\
7	0.563872705858414\\
8	0.559432381277899\\
9	0.550955977434403\\
10	0.542675416804365\\
11	0.551401790355387\\
12	0.54915980808043\\
13	0.537262898611903\\
14	0.546669488720833\\
15	0.557150469853094\\
16	0.563731980589134\\
17	0.561606952528243\\
18	0.554240137577236\\
19	0.527501396432295\\
19	0.527501396432295\\
18	0.554240137577236\\
17	0.561606952528243\\
16	0.563731980589134\\
15	0.557150469853094\\
14	0.546669488720833\\
13	0.537262898611903\\
12	0.54915980808043\\
11	0.551401790355387\\
10	0.542675416804365\\
9	0.550955977434403\\
8	0.559432381277899\\
7	0.563872705858414\\
6	0.565560852285104\\
5	0.553250337440242\\
4	0.563567207208195\\
3	0.584566456075738\\
2	0.64100994859013\\
1	0.738591230987051\\
}--cycle;

\addplot[area legend, draw=black, fill=black, fill opacity=0.2, forget plot]
table[row sep=crcr] {%
x	y\\
1	0.738591230987051\\
2	0.642112454619908\\
3	0.582558830035727\\
4	0.558526896543968\\
5	0.546393596811143\\
6	0.556757482334975\\
7	0.555002456237832\\
8	0.550671614211642\\
9	0.542429379679178\\
10	0.534219092369346\\
11	0.542083183957094\\
12	0.540073660176655\\
13	0.528750629421901\\
14	0.537310322087951\\
15	0.547380014241206\\
16	0.553985910149654\\
17	0.55234943282698\\
18	0.545452425923719\\
19	0.519883315557908\\
19	0.612536647258844\\
18	0.63608742479616\\
17	0.642386541028871\\
16	0.643802316473428\\
15	0.63764977484573\\
14	0.628340427514512\\
13	0.620549873532067\\
12	0.630987943398583\\
11	0.632805708756496\\
10	0.625645518520773\\
9	0.633228873171658\\
8	0.640689923105775\\
7	0.644382127056871\\
6	0.645309120870074\\
5	0.634257867592339\\
4	0.641690498312399\\
3	0.654432223312138\\
2	0.684317304955061\\
1	0.738591230987051\\
}--cycle;
\addplot [color=red, mark size=2.0pt, mark=asterisk, mark options={solid, red}, forget plot]
  table[row sep=crcr]{%
1	0.738591230987051\\
2	0.670011106792998\\
3	0.634655102485354\\
4	0.621639969365157\\
5	0.614570598562949\\
6	0.627910014774406\\
7	0.626530072730356\\
8	0.622320173869437\\
9	0.614156752886255\\
10	0.606208201572478\\
11	0.614934458322665\\
12	0.612674267656527\\
13	0.601031755142068\\
14	0.610405125566673\\
15	0.620613522080763\\
16	0.626882448112976\\
17	0.624636163799413\\
18	0.617409064835371\\
19	0.591338763779114\\
};
\addplot [color=black, mark size=1.0pt, mark=square*, mark options={solid, fill=black, black}, forget plot]
  table[row sep=crcr]{%
1	0.738591230987051\\
2	0.728841166842235\\
3	0.718284762909166\\
4	0.708805049005058\\
5	0.699930789205317\\
6	0.693856584847881\\
7	0.687653596116396\\
8	0.681691753509237\\
9	0.675085591968924\\
10	0.668318799459528\\
11	0.663511989930502\\
12	0.658259412452537\\
13	0.652116020177191\\
14	0.648720409986208\\
15	0.646517098640286\\
16	0.644886358270061\\
17	0.64259184455838\\
18	0.639572766635035\\
19	0.633771992658062\\
};
\end{axis}
\end{tikzpicture}
     \end{subfigure}
     \begin{subfigure}{0.45\textwidth}
         \centering
%
%
\begin{tikzpicture}
\pgfplotsset{every tick label/.append style={font=\small}}

\begin{axis}[%
width=6.5cm,
height=2.2cm,
at={(3.327in,0.481in)},
scale only axis,
xmin=1,
xmax=18,
xtick={3,7,11,15},
xticklabels={{1970},{1974},{1978},{1982}},
xlabel style={font=\bfseries\color{white!15!black}},
xlabel={Years},
ymin=0.6,
ymax=0.8,
ylabel style={font=\bfseries\color{white!15!black},shift={(0.0,-0.05)}},
ylabel={Pepsi's Beliefs},
axis background/.style={fill=white},
xmajorgrids,
ymajorgrids
]
\addplot [color=blue, mark size=1.0pt, mark=*, mark options={solid, fill=blue, blue}, forget plot]
  table[row sep=crcr]{%
1	0.738591230987051\\
2	0.744620276849233\\
3	0.743087144080133\\
4	0.744480014058034\\
5	0.734562817374239\\
6	0.718587384502697\\
7	0.719016091196367\\
8	0.719341418235089\\
9	0.723257351397837\\
10	0.657583595317576\\
11	0.660627889831085\\
12	0.651520953831866\\
13	0.656345693598383\\
14	0.628155103863401\\
15	0.608293382967169\\
16	0.614894888824271\\
17	0.632923740112681\\
18	0.637850059562571\\
19	0.624847444924705\\
};
\addplot [color=gray, mark size=2.0pt, mark=asterisk, mark options={solid, gray}, forget plot]
  table[row sep=crcr]{%
1	0.738591230987051\\
2	0.723164797799251\\
3	0.7125884120376\\
4	0.703818192602956\\
5	0.69554222974023\\
6	0.697483153123564\\
7	0.695672496425388\\
8	0.692228397549027\\
9	0.685775919414763\\
10	0.677826170757389\\
11	0.677589901909966\\
12	0.674419315020771\\
13	0.665943540511499\\
14	0.666806404221368\\
15	0.671100537898045\\
16	0.675800382938792\\
17	0.677174236367933\\
18	0.674827081834663\\
19	0.660077331708371\\
};
\addplot [color=gray, mark size=2.0pt, mark=asterisk, mark options={solid, gray}, forget plot]
  table[row sep=crcr]{%
1	0.738591230987051\\
2	0.736937638344537\\
3	0.742521002759345\\
4	0.746160152423113\\
5	0.747667017080943\\
6	0.757757124548423\\
7	0.761945790225164\\
8	0.763036667098502\\
9	0.76043222883283\\
10	0.755738094097856\\
11	0.758464737174936\\
12	0.757636972966625\\
13	0.751005263742112\\
14	0.753520786658676\\
15	0.759088943009413\\
16	0.764597244264613\\
17	0.766458027626876\\
18	0.764440040933402\\
19	0.750011467339333\\
};
\addplot [color=gray, mark size=2.0pt, mark=asterisk, mark options={solid, gray}, forget plot]
  table[row sep=crcr]{%
1	0.738591230987051\\
2	0.733831217214344\\
3	0.737411277362671\\
4	0.739317918174566\\
5	0.739189293766818\\
6	0.749619316405279\\
7	0.753094461562949\\
8	0.753163036853888\\
9	0.749063918682133\\
10	0.742836662132996\\
11	0.745750694964878\\
12	0.744412965109423\\
13	0.736322903783623\\
14	0.73940235696028\\
15	0.74599106761304\\
16	0.752316971752187\\
17	0.75414289543198\\
18	0.751383512855946\\
19	0.734076634169156\\
};
\addplot [color=red, mark size=2.0pt, mark=asterisk, mark options={solid, red}, forget plot]
  table[row sep=crcr]{%
1	0.738591230987051\\
2	0.726626475191061\\
3	0.717440231018757\\
4	0.709417818631144\\
5	0.701673881516855\\
6	0.701281096400047\\
7	0.698524229122589\\
8	0.694754503247986\\
9	0.688772943526957\\
10	0.681619819890792\\
11	0.679933288545543\\
12	0.676268427452149\\
13	0.668833063647978\\
14	0.668168427457747\\
15	0.670213227220454\\
16	0.672874594536409\\
17	0.673329546706514\\
18	0.671179075630164\\
19	0.659967939874302\\
};

\addplot[area legend, draw=black, fill=black, fill opacity=0.2, forget plot]
table[row sep=crcr] {%
x	y\\
1	0.738591230987051\\
2	0.723164797799251\\
3	0.7125884120376\\
4	0.703818192602956\\
5	0.69554222974023\\
6	0.697483153123564\\
7	0.695672496425388\\
8	0.692228397549027\\
9	0.685775919414763\\
10	0.677826170757389\\
11	0.677589901909966\\
12	0.674419315020771\\
13	0.665943540511499\\
14	0.666806404221368\\
15	0.671100537898045\\
16	0.675800382938792\\
17	0.677174236367933\\
18	0.674827081834663\\
19	0.660077331708371\\
19	0.750011467339333\\
18	0.764440040933402\\
17	0.766458027626876\\
16	0.764597244264613\\
15	0.759088943009413\\
14	0.753520786658676\\
13	0.751005263742112\\
12	0.757636972966625\\
11	0.758464737174936\\
10	0.755738094097856\\
9	0.76043222883283\\
8	0.763036667098502\\
7	0.761945790225164\\
6	0.757757124548423\\
5	0.747667017080943\\
4	0.746160152423113\\
3	0.742521002759345\\
2	0.736937638344537\\
1	0.738591230987051\\
}--cycle;

\addplot[area legend, draw=black, fill=black, fill opacity=0.2, forget plot]
table[row sep=crcr] {%
x	y\\
1	0.738591230987051\\
2	0.736937638344537\\
3	0.742521002759345\\
4	0.746160152423113\\
5	0.747667017080943\\
6	0.757757124548423\\
7	0.761945790225164\\
8	0.763036667098502\\
9	0.76043222883283\\
10	0.755738094097856\\
11	0.758464737174936\\
12	0.757636972966625\\
13	0.751005263742112\\
14	0.753520786658676\\
15	0.759088943009413\\
16	0.764597244264613\\
17	0.766458027626876\\
18	0.764440040933402\\
19	0.750011467339333\\
19	0.750011467339333\\
18	0.764440040933402\\
17	0.766458027626876\\
16	0.764597244264613\\
15	0.759088943009413\\
14	0.753520786658676\\
13	0.751005263742112\\
12	0.757636972966625\\
11	0.758464737174936\\
10	0.755738094097856\\
9	0.76043222883283\\
8	0.763036667098502\\
7	0.761945790225164\\
6	0.757757124548423\\
5	0.747667017080943\\
4	0.746160152423113\\
3	0.742521002759345\\
2	0.736937638344537\\
1	0.738591230987051\\
}--cycle;
\addplot [color=black, mark size=1.0pt, mark=square*, mark options={solid, fill=black, black}, forget plot]
  table[row sep=crcr]{%
1	0.738591230987051\\
2	0.727997095634584\\
3	0.717960536205325\\
4	0.708833477642199\\
5	0.700004014722245\\
6	0.695867091100962\\
7	0.690654315923835\\
8	0.685118626931026\\
9	0.678158093759284\\
10	0.670478538165937\\
11	0.666191002319796\\
12	0.660768442977824\\
13	0.653156901424877\\
14	0.650060810950535\\
15	0.648897671862209\\
16	0.648464277033912\\
17	0.64674099644543\\
18	0.64349079895333\\
19	0.634663284657312\\
};
\end{axis}
\end{tikzpicture}%
     \end{subfigure}
     \caption{Coke’s and Pepsi’s beliefs on market share evolution of Coke. The true market share evolution (blue line), the OLS regression fit (black line), the fit from the $L_2$-based variant (red line) and the fit from the polyhedral solution set (grey region). 
     }
  \label{fig:market_share_evolution}
\end{figure}
%
 On the other hand, the objective of the OLS regression is to find the parameters that lead to the best fit of the market share evolution, shown in Figure~\ref{fig:market_share_evolution}\footnote{We added parameter estimates that are at most 5\% more irrational than $\hat{\epsilon}_i$ to the polyhedral solution set (in this case a singleton).}. 
 In this case, the estimates given by our method are less in agreement with the actual market share evolution. This might be due to the fact that the market share is not only affected by the advertising expenditure but also by other economic factors and that the market share evolution is more complex than captured by the simplified Lanchester (resp. Sorger) model. Despite differences in the estimates (whose accuracy can only be computed by means of ground truth values of the parameters) and in their performances under different metrics, the above example has shown the applicability of our method to a real-world case scenario.


%

\section{Conclusion}

We designed a data-driven inverse-game-theoretic inference method to learn parametric utilities in multi-agent games.
Our proposed technique works with both Nash equilibrium profiles or better-response trajectories observations, is computationally tractable, and produces polyhedral solution sets containing all the parameters which rationalize the observed behaviour best.
%
The proposed dynamic inference method has been used to estimate Coke’s and Pepsi’s beliefs on market share evolution from their historical advertising expenditures. Unlike other standard inference methods the market share beliefs estimated with our approach minimize irrationality of the firms' observed expenditures.

\bibliography{ifacconf}             

\appendix

\section{Appendix}\label{sec:appendix}

\subsection{Proof of Proposition \ref{prop:1}}
\label{ap:proofOfProp1}
By definition of the irrationality loss function in \eqref{distance_func}, $\bar \epsilon:=\max_{i\in\mathcal{I}} \{d_i(\hat{\theta}_i; \mc D_i)\}$ implies that $0\leq d_i(\hat{\theta}_i; \mc D_i)\leq\bar\epsilon$, for all $i\in\mathcal{I}$.
Furthermore, by definition of the error function in \eqref{pre_error_func}, $0\leq d_i(\hat{\theta}_i;\mc D_i)\leq\bar\epsilon$ implies that $U_i(\bar{a}_i^j,\bar{a}_{-i}^j;\hat{\theta}_i)-U_i(a_i^j,\bar{a}_{-i}^j;\hat{\theta}_i)\leq\bar\epsilon$, $ \forall j\in\left[m\right]$.
\begin{enumerate}[label=(\roman*)]
    \item Under Assumption~\ref{ass:BRNEcond}~(i), the former inequality reads as $U_i(a_i(j),a_{-i}(j);\hat{\theta}_i)- U_i(a_i(j+1),a_{-i}(j);\hat{\theta}_i)\leq\bar\epsilon$,  $\forall j\in\left[m\right]$ and $\forall i\in\mathcal{I}$, which corresponds to the better-response condition in \eqref{better-response}, $\forall j\in\left[m\right]$ and $\forall i\in\mathcal{I}$. We conclude that $(a_i(j))_{j\in\left[m\right]}$ is an $\bar\epsilon$-better-response dynamics for each agent $i\in\mathcal{I}$.
    \item Under Assumption~\ref{ass:BRNEcond}~(ii), and the additional assumption $\bigcup_{j=1}^m \tilde{a}_i^j=\mathcal{A}_i$, the former inequality reads as $U_i(\tilde{a}_i^j,a_{-i}^*;\hat{\theta}_i)- U_i(a_i^*,a_{-i}^*;\hat{\theta}_i)\leq\bar\epsilon$,  $\forall \tilde{a}_i^j\in\mathcal{A}_i$, which corresponds to the best-response condition in \eqref{best-response}.
    By collecting the latter inequalities over all agents $i\in \mathcal{I}$, we conclude that $a^*$ is an $\bar\epsilon$-NE.
\end{enumerate}
{\hfill $\blacksquare$}


\subsection{Proof of Proposition \ref{LP_reformulation}}
\label{ap:ProofOfLPref}
Fix $\overline{\theta}_i\in\Theta_i $, where $\Theta_i \neq \varnothing$. Denote by $\Tilde{\epsilon}_i$  the maximum of $e_i(a_i^{j,k},\bar{a}^{j,k};\theta_i,\xi^k)$ over all data points $\{ \mc D_i^k \}_{k \in [n]}$ and set $\overline{\epsilon}_i=\max\{0,\Tilde{\epsilon}_i\}$. Thus, $\overline{\theta}_i\in\Theta_i$ and $\overline{\epsilon}_i\geq 0$ satisfy $e_i(a_i^{j,k},\bar{a}^{j,k};\theta_i,\xi^k)\leq\Tilde{\epsilon}_i\leq\overline{\epsilon}_i$ for all $j\in\left[m^k\right]$ and $k\in\left[n\right]$. Then, it follows by \S 4.1 in \citep{boyd04} that the LP \eqref{eq:pr3} is feasible.
{\hfill $\blacksquare$}

\end{document}